\documentclass[12pt]{amsart}

\usepackage{amssymb}

\theoremstyle{plain}
\newtheorem{thm}{Theorem}[section]
\newtheorem{lem}{Lemma}[section]
\newtheorem{cor}{Corollary}[section]
\newtheorem{prop}{Proposition}[section]

\theoremstyle{definition}

\numberwithin{equation}{section}

\begin{document}
\baselineskip=17pt
\title{On the Radius in Cayley--Dickson Algebras}

\author{Moshe Goldberg}

\address{Department of Mathematics,
Technion -- Israel Institute of Technology,
Haifa 32000, Israel}

\email{mg@technion.ac.il}

\author{Thomas J. Laffey}

\address{School of Mathematical Sciences,
University College Dublin,
Dublin 4, Ireland}

\email{thomas.laffey@ucd.ie}



\begin{abstract}
In the first two sections of this paper we provide a brief account of the Cayley--Dickson
algebras and prove that the {\em radius} on these algebras is given by the Euclidean norm.
With this observation we resort to three related topics: a variant of the Gelfand formula,
stability of subnorms, and the functional power equation.
\end{abstract}
\maketitle

\section{\label{sec:1}The Cayley--Dickson algebras}

The Cayley--Dickson algebras constitute a familiar series of algebras,
$\mathcal{A}_0, \mathcal{A}_1, \mathcal{A}_2,\ldots$ over the reals, with $\dim\mathcal{A}_{n}=2^n$.

The first five algebras in this series are the reals $\mathbb{R}$, the complex numbers
$\mathbb{C}$, the quaternions $\mathbb{H}$, the octonions $\mathbb{O}$, and the sedenions
$\mathbb{S}$. While $\mathbb{R}$ and $\mathbb{C}$ are both commutative and associative,
$\mathbb{H}$ is no longer commutative, and $\mathbb{O}$ and $\mathbb{S}$ are not even
associative. As it is, $\mathbb{O}$ is alternative, and $\mathbb{S}$ is merely power-associative.

As usual, we call an algebra $\mathcal{A}$ {\em alternative} if the subalgebra generated by any
two elements is associative. Further, $\mathcal{A}$ is called {\em power-associative} if the
subalgebra generated by any one element is associative. It thus follows that
{\em if} $\mathcal{A}$ {\em is power-associative, then the powers of each element in}
$\mathcal{A}$ {\em are unambiguously defined}.

It is well known that the Cayley--Dickson algebras can be inductively obtained from each other
by the following {\em Cayley--Dickson doubling process} (e.g., [W1]). We initiate this process
by setting $\mathcal{A}_0=\mathbb{R}$ and defining $a^\ast$, the {\em conjugate} of a real number
$a$, to equal $a$. Then, assuming that $\mathcal {A}_{n-1}$, $n \geq 1$, has been determined,
we define $\mathcal{A}_n$ to be the set of all ordered pairs
$$
\mathcal{A}_n = \{(a,b): ~ a,b \in \mathcal{A}_{n-1}\},
$$
such that addition and scalar multiplication are taken componentwise on the Cartesian product
$\mathcal{A}_{n-1}\times\mathcal{A}_{n-1}$, conjugation is determined by
$$
(a,b)^\ast=(a^\ast,-b),
$$
and multiplication is given by
$$
(a,b)(c,d)=(ac-d^\ast b, da+bc^\ast).
$$

With this definition, each element in $\mathcal{A}_n$ is of the form
$a=(\alpha_1,\ldots,\alpha_{2^n})$, $\alpha_j\in\mathbb{R}$. Moreover, it readily follows
that the distributive laws hold, and that the conjugate of $a$ and the unit element in
$\mathcal{A}_n$ are given, respectively, by
$$
a^\ast=(\alpha_1, -\alpha_2, \ldots, -\alpha_{2^n}) \quad \textrm{and} \quad
\textbf{1}_n=(1, 0, \ldots, 0).
$$

For example, we have
$$
\mathcal{A}_1=\{(\alpha,\beta): ~ \alpha,\beta \in \mathbb{R}\},
$$
where
$$
(\alpha,\beta)^\ast=(\alpha^\ast,-\beta), \quad
(\alpha,\beta)(\gamma,\delta)=(\alpha\gamma-\delta\beta,\delta\alpha+\beta\gamma), \quad
\textbf{1}_1=(1,0);
$$
hence we may identify $\mathcal{A}_1$ with $\mathbb{C}$ upon writing $z=\alpha+i\beta$ as
$(\alpha,\beta)$.

By construction, $\mathcal{A}_{n-1}$ {\em can be viewed as a subalgebra of} $\mathcal{A}_n$.
Since the octonions are non-associative, it follows that $\mathcal{A}_n$ is non-associative for
all $n\geq 3$. Similarly, since the sedenions are not alternative, $\mathcal{A}_n$ is not alternative for $n\geq 4$. It is known, however, that {\em all Cayley--Dickson algebras are power-associative}
---a fact that will be revisited in Theorem 1.1 below.

Carrying on with our short account, we shall now state the following known result whose proof is
provided for the reader's convenience:
\begin{lem}
\label{lem:1.1}
Each $a=(\alpha_1,\alpha_2,\ldots,\alpha_{2^n})\in\mathcal{A}_n$ is annihilated by the quadratic
polynomial
$$
p_a(t)=t^2-2\alpha_1t+|a|^2,
$$
where
$$
|a|=\sqrt{\alpha_1^2 + \cdots + \alpha_{2^n}^2}
$$
is the Euclidean norm.
\end{lem}

\noindent {\em Proof.} Let us show by induction that
\begin{equation}
\label{eq:1.1}
a^\ast a = |a|^2\textbf{1}_n, \quad a\in\mathcal{A}_n.
\end{equation}
The case $n=0$ is trivial, so assuming the assertion for $n-1$ and writing $a\in\mathcal{A}_n$ as
$$
a=(b,c), \quad b,c\in\mathcal{A}_{n-1},
$$
we get
\begin{equation*}
\begin{split}
a^\ast a & = (b,c)^\ast (b,c) = (b^\ast, -c)(b,c)
= (b^\ast b + c^\ast c, cb^\ast - cb^\ast) \\
& = (b^\ast b + c^\ast c, 0)
= (|b|^2 + |c|^2)(\textbf{1}_{n-1}, 0) = |a|^2 \textbf{1}_n.
\end{split}
\end{equation*}
We also have,
\begin{equation}
\label{eq:1.2}
a + a^\ast=(\alpha_1, \alpha_2, \ldots, \alpha_{2^n})
+ (\alpha_1, -\alpha_2, \ldots, -\alpha_{2^n})=(2\alpha_1, 0, \ldots, 0)=2\alpha_1\textbf{1}_n.
\end{equation}
So by (1.1) and (1.2),
$$
a^2 - 2\alpha_1 a + |a|^2\textbf{1}_n = a^2 - (a + a^\ast)a + a^\ast a =0,
$$
and we are done.
\qed

\medskip

Aided by this lemma one can prove that $\mathcal{A}_n$ is power-associative. Since this property
lies at the heart of our paper, and since the proof is short and not readily available in the
literature, we take the liberty of posting it here:

\begin{thm}
\label{thm:1.1} The Cayley--Dickson algebras are power-associative.
\end{thm}

\noindent {\em Proof.} By Lemma 3 in [A], $\mathcal{A}_n$ is power-associative if and only if
$$
a^2a=aa^2 \quad \textrm{and} \quad a^2a^2=(a^2a)a \quad \textrm{for all } a\in\mathcal{A}_n.
$$
By Lemma 1.1, each $a\in\mathcal{A}_n$ satisfies a relation of the form
\begin{equation}
a^2=\alpha a + \beta \textbf{1}_n \quad \alpha,\beta\in\mathbb{R}.
\end{equation}
\label{eq:1.3}
So,
\begin{equation}
\label{eq:1.4}
a^2a=(\alpha a + \beta \textbf{1}_n)a = \alpha a^2 + \beta a
= a(\alpha a + \beta \textbf{1}_n) = a a^2.
\end{equation}
Furthermore, by (1.3) and (1.4),
\begin{equation*}
\begin{split}
a^2 a^2 &= (\alpha a + \beta \textbf{1}_n)a^2
= \alpha a a^2 + \beta a^2 = \alpha a^2 a + \beta a^2 \\
&= (\alpha a^2 + \beta a)a = ((\alpha a + \beta \textbf{1}_n)a)a = (a^2a)a,
\end{split}
\end{equation*}
which yields the desired result.
\qed

\medskip

We mention, in passing, that it is a simple matter to verify that the sedenions are not alternative. Indeed, denoting by $e_j$, $j=1,\ldots,16$, the basis element of $\mathbb{S}$
whose $j$th entry is 1 and all others are zero, we put $o_1=e_4+e_{11}$ and $o_2=e_7-e_{16}$. Hence, consulting the multiplication table for the sedenions (e.g., [W2]), we find that
$o_1(o_1 o_2)=0$ while $o_1^2 o_2 = -2o_2$; so the option of alternativity is shattered.

\section{\label{sec:2}The radius on the Cayley--Dickson algebras}

We begin this section by considering an arbitrary finite-dimensional power-associative algebra
$\mathcal{A}$ over a field $\mathbb{F}$.

As is customary, by a {\em minimal polynomial} of an element $a$ in $\mathcal{A}$, we mean a monic
polynomial of lowest positive degree with coefficients in $\mathbb{F}$ that annihilates $a$.

With this familiar definition, we recall:
\begin{thm}
[{[G1, Theorem 1.1]}]
\label{thm:2.1}
Let $\mathcal{A}$ be a finite-dimensional power-associative algebra over a field $\mathbb{F}$. Then:

{\em (a)} Every element $a\in \mathcal{A}$ possesses a unique minimal polynomial.

{\em (b)} The minimal polynomial of $a$ divides every other polynomial that annihilates $a$.
\end{thm}

Restricting attention to the case where the base field of our algebra is either $\mathbb{R}$ or
$\mathbb{C}$, and denoting the minimal polynomial of an element $a\in \mathcal{A}$ by $q_a$,
we recall, [G1], that the {\em radius} of $a$  is defined to be the nonnegative quantity
$$
r(a)=\max \{ |\lambda|:~\lambda \in \mathbb{C}, \lambda \textrm{ is a root of } q_a \}.
$$

The radius has been computed for elements in several well-known finite-dimensional power-associative
algebras. For instance, it was shown in [G1, Section 2 and page 4072] that
the radius on the low dimensional Cayley--Dickson algebras $\mathbb{C}$, $\mathbb{H}$,
and $\mathbb{O}$ is given, in each case, by the corresponding Euclidean norm.

Another example of considerable interest emerged in [G1], where it was observed that
if $\mathcal{A}$ is an arbitrary matrix algebra over $\mathbb{R}$ or $\mathbb{C}$
(with the usual matrix operations), then the radius of a matrix $A \in \mathcal{A}$ is given
by the classical spectral radius,
$$
\rho(A)=\max \{ |\lambda|:~\lambda \in \mathbb{C}, \lambda \textrm{ is an eigenvalue of } A \}.
$$

With this last example in mind, we register the following theorem which tells us that the radius
retains some of the most basic properties of the spectral radius not only for matrices,
but in the general case as well.

\begin{thm}
[{[G1, Theorems 2.1 and 2.4]}]
\label{thm:2.2}
Let $\mathcal{A}$ be a finite-dimensional power-associative algebra over a field $\mathbb{F}$,
either $\mathbb{R}$ or $\mathbb{C}$. Then:

{\em (a)} The radius $r$ is a nonnegative function on $\mathcal{A}$.

{\em (b)} The radius is homogeneous, i.e., for all $a\in \mathcal{A}$ and $\alpha \in \mathbb{F}$,
$$
r(\alpha a)=|\alpha| r(a).
$$

{\em (c)} For all $a\in \mathcal{A}$ and all positive integers $k$,
$$
r(a^k)=r(a)^k.
$$

{\em (d)} The radius vanishes only on nilpotent elements of $\mathcal{A}$.

{\em (e)} The radius is a continuous function on $\mathcal{A}$.\footnote{Naturally,
a real-valued function on a finite-dimensional algebra $\mathcal{A}$ is said to be
{\em continuous} if it is continuous with respect to the unique finite-dimensional
topology on $\mathcal{A}$.}
\end{thm}

Having stated this theorem which shows that the radius is a natural extension of the spectral
radius, we are now ready to readdress the Cayley--Dickson algebras.

\begin{thm}
\label{thm:2.3}
The radius $r$ of an element $a = (\alpha_1,\ldots,\alpha_{2^n}) \in \mathcal{A}_n$ is the
Euclidean norm, i.e.,
$$
r(a)=|a|, \quad a\in \mathcal{A}_n.
$$
\end{thm}

\noindent {\em Proof.} Select $a$ in $\mathcal{A}_n$. By Lemma 1.1, the second order monic
$p_a(t) = t^2 - 2\alpha_1 t + |a|^2$ annihilates $a$. Hence $p_a$ is divisible by the minimal
polynomial of $a$. Further, since the roots of $p_a$ are
$$
t_\pm=\alpha_1 \pm \sqrt{\alpha_1^2 - |a|^2}
=\alpha_1 \pm \sqrt{-\alpha_2^2 - \cdots -\alpha_{2^n}^2},
$$
it follows that $|t_\pm|=|a|$; so $r(a)=|a|$ and the proof is complete.
\qed

\medskip

Since the radius on the Cayley--Dickson algebras is a norm, it automatically satisfies properties
(a), (b) and (e) of Theorem 2.2. Moreover, since a norm vanishes only on the zero element,
Theorem 2.2(d) implies that {\em the Cayley--Dickson algebras are void of nonzero nilpotent
elements.}

Combining Theorems 2.2(c) and 2.3, we realize that
$$
|a^k|=|a|^k \quad \textrm{ for all } a \in \mathcal{A}_n \textrm{ and } k=1,2,3,\ldots
$$
So if $a \neq 0$, then
$$
|a^k|=|a|^k \neq 0,
$$
which is another way of showing that the Cayley--Dickson algebras are void of nonzero nilpotents.

\section{\label{sec:3}Subnorms and a formula for the radius}

Again, let $\mathcal{A}$ be a finite-dimensional power-associative algebra over a field
$\mathbb{F}$, either $\mathbb{R}$ or $\mathbb{C}$. Following [GL2], we call a real-valued function
$$
f:\mathcal{A} \rightarrow \mathbb{R}
$$
a {\em subnorm} if for all $a \in \mathcal{A}$ and $\alpha \in \mathbb{F}$,
\begin{equation*}
\begin{split}
& f(a)>0, \quad a \neq 0, \\
& f(\alpha a)=|\alpha|f(a).
\end{split}
\end{equation*}

We recall that a real-valued function $N$ is a {\em norm} on $\mathcal{A}$ if for all
$a,b \in \mathcal{A}$ and $\alpha \in \mathbb{F}$,
\begin{equation*}
\begin{split}
& N(a)>0, \quad a \neq 0, \\
& N(\alpha a)=|\alpha|N(a), \\
& N(a+b)\leq N(a) + N(b);
\end{split}
\end{equation*}
hence, {\em a norm is a subadditive subnorm.} We also recall that in our finite-dimensional
settings, a norm is a continuous function on $\mathcal{A}$.

With the above definition of a subnorm, we may now state a variant of a result which is well
known in the context of complex Banach algebras (e.g., [R, Theorem 18.9], [L, Chapter 17,
Theorem 4]), and which is often referred to as the Gelfand formula.

\begin{thm}
[{[G2, Theorem 2.1]}]
\label{thm:3.1}
Let $f$ be a continuous subnorm on a a finite-dimensional power-associative algebra
$\mathcal{A}$ over $\mathbb{R}$ or $\mathbb{C}$. Then:
\begin{equation}
\label{eq:3.1}
\lim_{k\to\infty}f(a^k)^{1/k} = r(a) \quad \textit{for all } a\in \mathcal{A}.
\end{equation}
\end{thm}

In particular, for the Cayley--Dickson algebras we get:
\begin{cor}
\label{cor:3.1}
If $f$ is a continuous subnorm on $\mathcal{A}_n$, then
\begin{equation}
\label{eq:3.2}
\lim_{k\to\infty}f(a^k)^{1/k} = |a| \quad \textit{for all } a\in \mathcal{A}_n.
\end{equation}
\end{cor}

For example, we observe that for each fixed $p$, $0<p \leq \infty$, the function
\begin{equation}
\label{eq:3.3}
|a|_p=(|\alpha_1|^p + \cdots + |\alpha_{2^n}|^p)^{1/p},\quad
a=(\alpha_1,\ldots,\alpha_{2^n})\in \mathcal{A}_n,
\end{equation}
is a continuous subnorm on $\mathcal{A}_n$ (a norm precisely when $1 \leq\ p \leq \infty$).
Thus, we get
$$
\lim_{k\to\infty}|a^k|_p^{1/k}=|a| \quad \textrm{for all } a\in \mathcal{A}_n.
$$

We point out that contrary to norms, a subnorm on an algebra $\mathcal{A}$ over $\mathbb{R}$ or $\mathbb{C}$ with $\dim \mathcal{A} \geq 2$, may fail to be continuous. For instance,
[G2, Section 3], let $f$ be a continuous subnorm on such an algebra. Fix an element
$a_0 \neq 0$ in $\mathcal{A}$, and consider
$$
\textbf{V}=\left\{ \alpha a_0: ~\alpha \in \mathbb{F} \right\},
$$
the linear subspace of $\mathcal{A}$ generated by $a_0$. For each real $\kappa$, $\kappa >1$,
define
\begin{equation}
\label{eq:3.4}
g_\kappa(a)=
\left
\{\begin{array}{ll}
\kappa f(a), & a \in \textbf{V}, \\
f(a)       , & a \in \mathcal{A} \smallsetminus \textbf{V}.
\end{array}\right.
\end{equation}
Then evidently, $g_\kappa$ is a subnorm on $\mathcal{A}$. Further,
$g_\kappa$ is discontinuous at $a_0$ since
$$
\lim_{\substack{a \to a_0 \\ a \notin \mathbf{V}}} g_\kappa(a)
=\lim_{\substack{a \to a_0 \\ a \notin \mathbf{V}}} f(a)
=\lim_{a \to a_0} f(a)=f(a_0) \neq g_\kappa(a_0).
$$
Motivated by the above example, we shall next show that {\em a subnorm on a finite-dimensional
power-associative algebra over} $\mathbb{R}$ {\em or} $\mathbb{C}$ {\em may satisfy formula}
(3.1) {\em without being continuous}.

To this end, we recall that two subnorms $f$ and $g$ are {\em equivalent} on an algebra
$\mathcal{A}$ if there exist constants $\mu > 0$, $\nu >0$, such that
$$
\mu f(a) \leq g(a) \leq \nu f(a) \quad \textrm{for all } a\in\mathcal{A}.
$$

With this familiar definition, we may now register:

\begin{prop}
[{[G2, Theorem 3.1]}]
\label{prop: 3.1}
Let $g$ be a subnorm on a finite-dimensional power-associative algebra $\mathcal{A}$ over
$\mathbb{R}$ or $\mathbb{C}$. If $g$ is equivalent to a continuous subnorm on $\mathcal{A}$,
then $g$ satisfies formula {\em (3.1)}.
\end{prop}

To illustrate Proposition 3.1, let $\mathcal{A}$ be a finite-dimensional power-associative
algebra over $\mathbb{R}$ or $\mathbb{C}$ with $\dim \mathcal{A}=2$, and let us fall back on
the discontinuous subnorm $g_\kappa$ in (3.4). Since $\kappa >1$, we have
\begin{equation}
\label{eq:3.5}
f(a) \leq g_\kappa(a) \leq \kappa f(a), \quad a \in \mathcal{A};
\end{equation}
so (as indicated in Section 3 of [G2]), our proposition yields
$$
\lim_{k\to\infty}g_\kappa(a^k)^{1/k}=r(a) \quad \textrm{for all } a\in \mathcal{A}.
$$

In contrast with the last example, we shall now show that {\em not all discontinuous
subnorms on Cayley--Dickson algebras satisfy formula} (3.2). In order to support this
assertion, consider the familiar Cauchy equation
\begin{equation}
\label{eq:3.6}
\varphi(x+y) = \varphi(x) + \varphi(y), \quad x,y\in \mathbb{R},
\end{equation}
whose (real) solutions have been discussed in the literature for over a century (e.g., [H],
[HLP, Section 3.20], [HR], [B, Section 20], and [GL3, Section 2]). It is well known that
every solution of (3.6) satisfies
\begin{equation}
\label{eq:3.7}
\varphi(\gamma x)=\gamma \varphi(x) \quad \textrm{for all rational } \gamma \textrm{ and real } x.
\end{equation}
Hence, the only continuous solutions of (3.6) are of the form
$$
\varphi(x)=x \varphi(1), \quad x \in \mathbb{R},
$$
where $\varphi(1)$ is an arbitrary real value. It is also known that equation (3.6) has
discontinuous solutions, and that all such solutions are
discontinuous everywhere and unbounded (both from below and above) on any interval in $\mathbb{R}$.
Moreover, given a positive number $c$, one may select a discontinuous solution $\varphi$ with
$\varphi(c)=0$. Consequently,
\begin{equation}
\label{eq:3.8}
\varphi(x + c)=\varphi(x) + \varphi(c)=\varphi(x), \quad x\in \mathbb{R},
\end{equation}
so $\varphi$ can be chosen to be $c$-periodic. In fact, (3.8) tells us that $c$ is a period
of $\varphi$ if and only if $\varphi(c) = 0$; thus if $c$ is a period then, by (3.7), so is
every rational multiple of $c$.

With these remarks, we may quote the following result.

\begin{thm}
[{[GL3, Theorems 2.1(a,c) and 2.2(a)]}]
\label{thm:3.2}
Let $f$ be a continuous subnorm on $\mathbb{C}$, and let $\varphi$ be a discontinuous
$\pi$-periodic solution of the Cauchy equation {\em (3.6)}. Then:

{\em (a)} The function
\begin{equation}
\label{eq:3.9}
g_\varphi(z) = f(z)e^{\varphi(\arg z)}, \quad z\in \mathbb{C},
\end{equation}
where $\arg z$ denotes the principal argument of $z$, $0 \leq \arg z < 2\pi$, and where
$\arg 0 = 0$, is a subnorm on $\mathbb{C}$.

{\em (b)} $g_\varphi$ is discontinuous everywhere in $\mathbb{C}$.

{\em (c)} For every $z \in \mathbb{C}$,
\begin{equation}
\label{eq:3.10}
\lim_{k \to \infty} g_\varphi(z^k)^{1/k}=|z|e^{\varphi(\arg z)}.
\end{equation}
\end{thm}

Now, as indicated in Section 3 of [G2], equation (3.10) implies that the utterly discontinuous
subnorm $g_\varphi$ satisfies formula (3.2) only for those values of $z$ for which
\begin{equation}
\label{eq:3.11}
\varphi(\arg z)=0.
\end{equation}
As $\varphi$ is unbounded on any subinterval of $[0,2\pi)$, we infer that the set of points where the subnorm $g_\varphi$ violates formula (3.2) is dense in $\mathbb{C}$.

We further note that since our $\varphi$ is $\pi$-periodic, the remarks preceding Theorem 3.2
imply that (3.11) holds whenever $\arg z$ is a rational multiple of $\pi$. Hence the set of
complex points where $g_\varphi$ does satisfy formula (3.2) is also dense in $\mathbb{C}$.

We close this section by appealing to Theorems 3.1(a,c) and 3.2(a) in [GL3], by which we acquire
the following analog of Theorem 3.2 for the quaternions.

\begin{thm}
\label{thm:3.3}
Let $f$ be a continuous subnorm on $\mathbb{H}$, and let $\varphi$ be a discontinuous
$\pi$-periodic solution of the Cauchy equation {\em (3.6)}. Then:

{\em (a)} The function
\begin{equation}
\label{eq:3.12}
h_\varphi(q) = f(q)e^{\varphi(\arg q)}, \quad q\in \mathbb{H},
\end{equation}
where
$$
\arg q = \arg (\alpha + i \sqrt{\beta^2 + \gamma^2 + \delta^2}), \quad
q = \alpha + i\beta + j\gamma + k\delta \in \mathbb{H},
$$
is a subnorm on $\mathbb{H}$.

{\em (b)} $h_\varphi$ is discontinuous everywhere in $\mathbb{H}$.

{\em (c)} For every $q \in \mathbb{H}$,
$$
\lim_{k \to \infty} h_\varphi(q^k)^{1/k}=|q|e^{\varphi(\arg q)}.
$$
\end{thm}

With this result at hand, we may echo the argument which followed Theorem 3.2 and observe that
the set on which $h_\varphi$ satisfies formula (3.2) and the one on which $h_\varphi$ violates
this formula are both dense in $\mathbb{H}$.

\section{\label{sec:4}Stability of subnorms and radial dominance}

In accordance with routine nomenclature for norms, we follow [GL2] and say that a subnorm $f$
on a finite-dimensional power-associative algebra $\mathcal{A}$ over $\mathbb{R}$ or $\mathbb{C}$
is {\em stable} if there exists a positive constant $\sigma$ such that
\begin{equation}
\label{eq:4.1}
f(a^k) \leq \sigma f(a)^k \quad \textrm{for all } a \in \mathcal{A} \textrm{ and } k=1,2,3,\ldots
\end{equation}
Furthermore, we say that $f$ is {\em radially dominant} on $\mathcal{A}$ if $f$ majorizes the radius, i.e.,
$$
f(a) \geq r(a), \quad a\in \mathcal{A}.
$$

With these notions of stability and radial dominance, we may now cite one of the main results
in [G1].

\begin{thm}
[{[G1, Theorem 3.4]}]
\label{thm:4.1}
Let $f$ be a continuous subnorm on a finite-dimensional power-associative algebra $\mathcal{A}$ over $\mathbb{R}$ or $\mathbb{C}$, and let $\mathcal{A}$ be void of nonzero nilpotents. Then $f$
is stable if and only if $f$ is radially dominant on $\mathcal{A}$.
\end{thm}

By Theorem 2.2, in the absence of nonzero nilpotents the radius is a continuous subnorm. Hence, Theorem 4.1 instantly implies:

\begin{cor}
\label{cor:4.1}
Let $\mathcal{A}$, a finite-dimensional power-associative algebra over $\mathbb{R}$ or $\mathbb{C}$, be void of nonzero nilpotents. Then the radius is the smallest of all stable continuous subnorms on $\mathcal{A}$.
\end{cor}

Another immediate consequence of Theorem 4.1 reads:
\begin{cor}
\label{cor:4.2}
Let $f$ be a subnorm on the Cayley--Dickson algebra $\mathcal{A}_n$. Then $f$ is stable if and
only if $f$ majorizes the Euclidean norm on $\mathcal{A}_n$.
\end{cor}

We illustrate the last corollary by noticing that the subnorm in (3.3) is stable if and
only if
$$
|a|_p \geq |a| \quad \textrm{ for all } a \in \mathcal{A}_n,
$$
which holds, of course, only when $0 < p \leq 2$.

With this example and with Theorem 3.1 in the bag, we see that {\em while continuity of a
subnorm is enough to force formula} (3.1), {\em it is not enough to force stability, not even
in the absence of nonzero nilpotents}.

To further illustrate Corollary 4.2, let $w=(u,v)$ be a fixed vector of two positive entries,
and consider the weighted sup norm on $\mathbb{C}$,
$$
\|z\|_{w,\infty}=\max \{ u|\alpha|, v|\beta| \}, \quad z = \alpha + i\beta \in \mathbb{C}.
$$
By the last corollary, $\| \cdot \|_{w,\infty}$ is stable on $\mathbb{C}$ if and only if
$$
\max \{ u|\alpha|, v|\beta| \} \geq \sqrt{\alpha^2 + \beta^2} \quad
\textrm{ for all } \alpha, \beta \in \mathbb{R};
$$
which, by Theorem 3.1 in [GL2], holds precisely when
$$
u^2v^2 \geq u^2 + v^2.
$$

A similar example is obtained by considering the weighted $\ell_1$-norm,
\begin{equation}
\label{eq:4.2}
\|z\|_{w,1} = u|\alpha| + v|\beta|, \quad z = \alpha + i\beta \in \mathbb{C},
\end{equation}
where again, $w=(u,v)$ is a fixed 2-vector of positive entries. By Corollary 4.2,
$\|\cdot\|_{w,1}$ is stable on $\mathbb{C}$ if and only if
$$
u|\alpha| + v|\beta| \geq \sqrt{\alpha^2 + \beta^2} \quad
\textrm{ for all } \alpha, \beta \in \mathbb{R};
$$
which, by Theorem 3.3 in [GL2], is valid exactly for
$$
u\geq1, v\geq1.
$$

To put Theorem 4.1 in perspective, we shall now register a simple result which implies that
{\em under certain conditions stability and radial dominance may prevail even when continuity
does not.}

\begin{prop}
[{[G3, Proposition 5.1(b)]}]
\label{prop: 4.1}
Let $g$ be a stable subnorm on a finite-dimensional power-associative algebra $\mathcal{A}$ over
$\mathbb{R}$ or $\mathbb{C}$. If $g$ is equivalent to a continuous subnorm on $\mathcal{A}$,
then $g$ is radially dominant.
\end{prop}

Indeed, let $\mathcal{A}$ be a finite-dimensional power-associative algebra over $\mathbb{R}$ or $\mathbb{C}$ with $\dim\mathcal{A} \geq 2$, let $f$ be a continuous stable subnorm
on $\mathcal{A}$, and let $g_\kappa$ be the discontinuous subnorm in (3.4). By (3.5) and (4.1),
$$
g_\kappa(a^k) \leq \kappa f(a^k) \leq \kappa\sigma f(a)^k \leq \kappa\sigma g(a)^k,\quad
a \in \mathcal{A},~k=1,2,3,\ldots
$$
So $g_\kappa$ is stable and, by Proposition 4.1, $g_\kappa$ is radially dominant a well.

In view of this example, we proceed to show that {\em not all stable subnorms on Cayley--Dickson algebras are radially dominant}. We begin by citing:

\begin{thm}
[{[GL3, Theorem 2.2(b)]}]
\label{thm:4.2}
Let $f$ be a continuous stable subnorm on $\mathbb{C}$, let $\varphi$ be a discontinuous $\pi$-periodic solution of equation {\em (3.6)}, and let $g_\varphi$ be the discontinuous
subnorm in {\em (3.9)}. Then $g_\varphi$ is stable on $\mathbb{C}$.
\end{thm}

Now, letting $f$, $\varphi$ and $g_\kappa$ be as in Theorem 4.2, and walking in the footsteps
of Section 4 in [G2], we fix $z_0 \in \mathbb{C}$, $z_0 \neq 0$, and select $\varepsilon$,
$0 < \varepsilon <|z_0|$. Since $f$ is continuous on $\mathbb{C}$, $f$ is bounded on the closed
disc $D_\varepsilon = \{z:~ |z-z_0|\leq\varepsilon\}$. Moreover, since $\varphi$ is unbounded
from below on $[0,2\pi)$, the exponential function $e^{\varphi(\arg z)}$ attains arbitrary small positive values in $D_\varepsilon(z_0)$; so there exists a point
$z_\varepsilon \in D_\varepsilon(z_0)$ such that
$$
g_\varphi(z_\varepsilon)=f(z_\varepsilon)e^{\varphi(\arg z_\varepsilon)} < |z_0| - \varepsilon.
$$
We thus obtain
\begin{equation}
\label{eq:4.3}
g_\varphi(z_\varepsilon) < |z_\varepsilon|,
\end{equation}
and it follows that the ubiquitously discontinuous stable subnorm $g_\varphi$ is not radially dominant on $\mathbb{C}$.

As a matter of fact, since $f(\alpha z) = |\alpha|f(z)$ for all $z \in \mathbb{C}$ and
$\alpha \in \mathbb{R}$, and since every rational multiple of $\pi$ ia a period of $\varphi$, it is not hard to conclude from (4.3) that the set of complex numbers where $g_\varphi$ fails to be radially dominant is dense in $\mathbb{C}$.

Furthermore, since $f$ vanishes only at $a=0$, $f$ is bounded away  from zero on
$D_\varepsilon(z_0)$; and since $\varphi$ is unbounded from above on open subintervals of
$[0,2\pi)$, we can find a point $w_\varepsilon \in D_\varepsilon(z_0)$ for which
$$
g_\varphi(w_\varepsilon)=f(w_\varepsilon)e^{\varphi(\arg w_\varepsilon)}
> |z_0| +\varepsilon > |w_\varepsilon|.
$$
Therefore, the set where the subnorm $g_\varphi$ majorizes the radius is also dense in $\mathbb{C}$, showing that no inequality between $g_\varphi$ and the radius is possible.

To obtain an analog of Theorem 4.2 for the quaternions, we resort to:

\begin{thm}
[{[GL3, Theorem 3.2(b)]}]
\label{thm:4.3}
Let $f$ be a continuous stable subnorm on $\mathbb{H}$, let $\varphi$ be a discontinuous
$\pi$-periodic solution of the Cauchy equation {\em (3.6)}, and let $h_\varphi$ be the discontinuous subnorm in {\em (3.12)}. Then $h_\varphi$ is stable on $\mathbb{H}$.
\end{thm}

Rephrasing the argument that followed Theorem 4.2, we can now gather that if $f$, $\varphi$ and $h_\varphi$ are as in Theorem 4.3, then $h_\varphi$ neither majorizes nor is majorized by the
radius on $\mathbb{H}$.

We end this section by recalling, [GL2], that a subnorm $f$ on a finite-dimensional power-associative algebra $\mathcal{A}$ over $\mathbb{R}$ or $\mathbb{C}$ is {\em strongly stable} if (4.1) holds with $\sigma=1,$ that is,
$$
f(a^k) \leq f(a)^k \quad \textrm{for all } a \in \mathcal{A} \textrm{ and } k=1,2,3,\ldots
$$

To illustrate this definition, we allude to Theorem 2.2 and observe that
{\em if $\mathcal{A}$, a finite-dimensional power-associative algebra over $\mathbb{R}$
or $\mathbb{C}$, is void of nilpotent elements, then the radius is a strongly stable subnorm
on $\mathcal{A}$. }

Another example is obtained by referring to Theorem 3.4 in [GL2], which tells us that the norm
$\| \cdot \|_{w,1}$ in (4.2) is strongly stable on $\mathbb{C}$ if and only if
$v^2 \geq u \geq 1$. Since we already know that $\| \cdot \|_{w,1}$ is stable only when $u \geq 1$ and $v \geq 1$, we deduce that $\| \cdot \|_{w,1}$ is stable but not strongly stable precisely for
$u \geq v^2 \geq 1$.

In light of this observation, it would be useful to obtain a characterization of strongly stable subnorms on the Cayley--Dickson algebras or, better yet, on general finite-dimensional power-associative algebras over $\mathbb{R}$ or $\mathbb{C}$.

\section{\label{sec:5}The power equation}

Let $\mathcal{A}$ be a finite-dimensional power-associative algebra over $\mathbb{R}$ or $\mathbb{C}$ and let
$$
f:\mathcal{A} \rightarrow \mathbb{R}
$$
be a real-valued function on $\mathcal{A}$. We say that $f$ is a {\em solution of the power equation} on $\mathcal{A}$ if
\begin{equation}
\label{eq:5.1}
f(a^k)=f(a)^k \quad \textrm{for all } a \in \mathcal{A} \textrm{ and } k=1,2,3,\ldots
\end{equation}
If, in addition, $f$ is a subnorm, we say that $f$ is a {\em submodulus} on $\mathcal{A}$.

Obviously, the only constant solutions of (5.1) are $f=0$ and $f=1$. When the base field of $\mathcal{A}$ is the real or the complex numbers, all we have to do in order to produce a non-trivial solution of (5.1), is to invoke Theorem 2.2 which yields:

\begin{cor}
\label{cor: 5.1}
Let $\mathcal{A}$ be a finite-dimensional power-associative algebra over $\mathbb{R}$ or $\mathbb{C}$. Then the radius $r$ is a nonnegative, homogeneous, and continuous solution of the power equation which vanishes only on nilpotent elements.
\end{cor}

In general, the radius is not a submodulus, nor even a subnorm, since it vanishes on nonzero nilpotents. Barring such nilpotents, we may obtain the following, more subtle version of
Corollary 5.1:

\begin{thm}
[{[G1, Theorem 3.3(b)]}]
\label{thm: 5.1}
Let $\mathcal{A}$, a finite-dimensional power-associative algebra over $\mathbb{R}$ or $\mathbb{C}$, be void of nonzero nilpotent elements. Then the radius is the only continuous submodulus on $\mathcal{A}$.
\end{thm}

Since the Cayley--Dickson algebras are void of nonzero nilpotents, this last theorem allows us to record:

\begin{cor}
\label{cor: 5.2}
The radius is the only continuous submodulus on $\mathcal{A}_n$.
\end{cor}

We remark that the definition of submodulus gives rise to yet another simple result which holds
for finite-dimensional as well as for infinite-dimensional algebras:

\begin{prop}
[{[GGL, Proposition 3]}]
\label{prop:5.1}
If $\mathcal{A}$, a power-associative algebra over $\mathbb{R}$ or $\mathbb{C}$, contains nonzero nilpotents, then $\mathcal{A}$ has no submodulus.
\end{prop}

With corollary 5.2 in mind, and recalling that the radius on the Cayley--Dickson algebras is the Euclidean norm, we appeal again to two of the results in [GL3] which tell us that solutions of
the power equation on $\mathbb{C}$ and $\mathbb{H}$ may fail to have even a shred of continuity:

\begin{thm}
[{[GL3, Theorems 2.1(b,c) and 3.1(b,c)]}]
\label{thm:5.2}
Let $\varphi$ be a discontinuous $\pi$-periodic solution of the Cauchy equation {\em (3.6)}. Then:

{\em (a)} The function
$$
g_\varphi(z)=|z|e^{\varphi(\arg z)}, \quad z\in \mathbb{C},
$$
is a submodulus on $\mathbb{C}$ which is discontinuous everywhere.

{\em (b)} Analogously, the function
$$
h_\varphi(z)=|q|e^{\varphi(\arg q)}, \quad q\in \mathbb{H},
$$
is a submodulus on $\mathbb{H}$ which is discontinuous everywhere.
\end{thm}

Whether continuous or not, we maintain that all solutions of the power equation on
$\mathcal{A}_n$, $n \geq 1$, are nonnegative. To establish this fact, we begin with the
following lemma which seems to be of independent interest.

\begin{lem}
\label{lem:5.1}
Every element in $\mathcal{A}_n$, $n \geq 1$, is a square; i.e., for each $a \in \mathcal{A}_n$ there is an element $b \in \mathcal{A}_n$ such that $a=b^2$.
\end{lem}

\noindent {\em Proof.} Select $a=(\alpha_1,0 \ldots, 0)$ in $\mathcal{A}_n$ with
$\alpha_1 \neq 0$. A simple induction on $n \geq 1$ implies that if $\alpha_1 \geq 0$, then $(\sqrt{\alpha_1}, 0, \ldots, 0)$ is a square root of $a$, and if $\alpha_1 < 0$, then
$(0, \sqrt{-\alpha_1}, 0, \ldots, 0)$ will do.

Next, select $a=(\alpha_1,\ldots,\alpha_{2^n})$ with $\alpha_j \neq 0$ for some $j \neq 1$.
Then $|a| > |\alpha_1|$; so $\alpha_1 + |a| > 0$, and
$$
b=\frac{1}{\sqrt{2\alpha_1 + 2|a|}}(a + |a|\textbf{1}_n)
$$
is well defined.

By Lemma 1.1,
$$
a^2=2\alpha_1 a - |a|^2 \textbf{1}_n.
$$
So,
\begin{equation*}
\begin{split}
b^2 &= \frac{1}{2\alpha_1 + 2|a|}(a + |a|\textbf{1}_n)^2
= \frac{1}{2\alpha_1 + 2|a|}(a^2 + 2|a|a + |a|^2\textbf{1}_n) \\
&= \frac{1}{2\alpha_1 + 2|a|}(2\alpha_1 a - |a|^2\textbf{1}_n
+ 2|a|a + |a|^2\textbf{1}_n) = a,
\end{split}
\end{equation*}
and the assertion is secured.
\qed

\medskip

We can now post:

\begin{thm}
[{Compare [GL1, Proposition 3.1]}]
\label{thm:5.3}
All solutions of the power equation on $\mathcal{A}_n$, $n \geq 1$, are nonnegative.
\end{thm}

\noindent {\em Proof.} Let $f$ be a solution of the power equation on $\mathcal{A}_n$,
$n \geq 1$. Select an element $a \in \mathcal{A}_n$. So by Lemma 5.1, $a = b^2$ for some
$b \in \mathcal{A}_n$. Whence,
$$
f(a)=f(b^2)=f(b)^2 \geq 0
$$
and the proposition follows.
\qed

\medskip

The function
$$
f(a)=a, \quad a \in \mathbb{R},
$$
shows, of course, that Theorem 5.3 is false for $\mathcal{A}_0$.

We comment that with the help of the solutions of the power equation which are already in our possession, it is easy to construct a plethora of others. For instance, if $f$ is a solution,
then so is $f^\kappa$ for any fixed positive constant $\kappa$. Also, if $f$ and $g$ are solutions, then so are $\max \{f,g \}$ and $\min \{f,g \}$, as well as $f^\kappa g^\tau$ for any pair of
fixed positive $\kappa$ and $\tau$.

We conclude this paper by recalling that on $\mathbb{R}$, $\mathbb{C}$, $\mathbb{H}$ and $\mathbb{O}$, the radius is {\em multiplicative}; i.e.,
\begin{equation}
\label{eq:5.2}
|ab|=|a| |b|, \quad a,b \in \mathcal{A}_n,~n=0,1,2,3,
\end{equation}
where for the octonions this identity follows from the Eight Square Theorem, [D]. This suffices
to imply, of course, that the Euclidean norm is a solution of the power equation on the reals,
the complex numbers, the quaternions, and the octonions. For $n \geq 4$, however, formula (5.2) fails since the Cayley--Dickson algebras contain zero divisors. To confirm this known fact, we
use once more, as we did at the end of Section 1, the octonions $o_1=e_4+e_{11}$ and $o_2=e_7-e_{16}$, and notice that $o_1 o_2=0$.


\begin{thebibliography}{-------}

\bibitem[A]{Albert}A.~A.~Albert,
{\em Power-associtive rings},
Trans. Amer. Math. Soc. $\mathbf{64}$ (1948), 552--593.

\bibitem[B]{Boas}R.~P.~Boas Jr.,
{\em A Primer of Real Functions},
Amer. Math. Soc., Providence, Rhode Island, 1960.

\bibitem[D]{Dickson}L.~E.~Dickson,
{\em On quaternions and their generalization and the history of the Eight Square Theorem},
Ann. of Math. $\mathbf{20}$ (1918--1919), 155--171.

\bibitem[G1]{Goldberg 1}M.~Goldberg,
{\em Minimal polynomials and radii of elements in finite-dimensional power-associative algebras},
Trans. Amer. Math. Soc. $\mathbf{359}$ (2007), 4055--4072.

\bibitem[G2]{Goldberg 2}M.~Goldberg,
\emph{Radii and subnorms on finite-dimensional power-associative algebras},
Linear Multilinear Algebra, $\mathbf{55}$ (2007), 405--415.

\bibitem[G3]{Goldberg 3}M.~Goldberg,
{\em Stable subnorms on finite-dimensional power-associative algebras},
Electron. J. Linear Algebra $\mathbf{17}$ (2008), 359--75.


\bibitem[GGL]{Goldberg-Guralnick-Luxemburg}M.~Goldberg,~R.~Guralnick and W.~A.~J.~Luxemburg,
{\em Stable subnorms II},
Linear Multilinear Algebra $\mathbf{51}$ (2003), 209--219.

\bibitem[GL1]{Goldberg-Levy}M.~Goldberg and E.~Levy,
{\em The power equation},
Electron. J. Linear Algebra $\mathbf{22}$ (2011), 810--821.

\bibitem[GL2]{Goldberg-Lux 1}M.~Goldberg and W.~A.~J.~Luxemburg,
{\em Stable subnorms},
Linear Algebra Appl. $\mathbf{307}$ (2000), 89--101.

\bibitem[GL3]{Goldberg-Lux 2}M.~Goldberg and W.~A.~J.~Luxemburg,
{\em Discontinuous subnorms},
Linear Multilinear Algebra $\mathbf{49}$ (2001), 1--24.

\bibitem[H]{Hamel}G.~Hamel,
{\em Eine Basis aller Zahlen und die unstetigen L\"{o}sungen der Functionalgleichung:} $f(x+y)=f(x)+f(y)$,
Math. Ann. $\mathbf{60}$ (1905), 450--462.

\bibitem[HLP]{Hardy-Littlewood-Polya}G.~H.~Hardy,~J.~E.~Littlewood and G.~P\'{o}lya,
{\em Inequalities}, Cambridge Univ. Press, Cambridge, 1934.

\bibitem[HR]{Hahn-Rosenthal}H.~Hahn and A.~Rosenthal,
{\em Set Functions},
Univ. of New Mexico Press, Albuquerque, 1948.

\bibitem[L]{Lax}P.~D.~Lax,
{\em Functional Analysis},
Wiley-Interscience, New York, 2002.

\bibitem[R]{Rudin}W.~Rudin,
{\em Real and Complex Analysis}, 3rd edition,
McGraw-Hill, New York, 1987.

\bibitem[W1]{Wikipedia 1}
http://en.wikipedia.org/wiki/Cayley-Dickson\_construction

\bibitem[W2]{Wikipedia 2}
http://en.wikipedia.org/wiki/Sedenions

\end{thebibliography}
\end{document}